\documentclass[a4paper,12pt]{article}   
  \usepackage{makeidx}
  \usepackage{graphicx,graphics,color}
\usepackage{latexsym,doc,fontenc,exscale}
\usepackage{amsmath,amsfonts,amsthm,eucal,amssymb}

\usepackage{multicol}
\usepackage{times}
\usepackage{pstricks}

\begin{document}

 \title{\vspace*{-2cm}\large Fourier transform for functions of bicomplex variables}%\vspace{2cm}
 %\vspace{17cm}
 %\noindent

\author{\small A Banerjee$^{1}$\footnote{abhijit.banerjee.81@gmail.com, abhijit\_banerjee@hotmail.com},
 S K Datta $^{2}$\footnote{sanjib\_kr\_datta@yahoo.co.in},
 Md. A Hoque $^{3}$\footnote{mhoque3@gmail.com}\\
 \small $^{1}$Department of Mathematics, Krishnath College,
                                        Berhampore, Murshidabad 742101, India\\
  \small $^{2}$ Department of Mathematics, University of Kalyani, Nadia 741235, India\\
 \small $^{3}$Gobargara High Madrasa (H.S.),Hariharpara,
                                    Murshidabad 742166,  India\\
 }
 \date{}
 %\today
 %\huge{\date{February 16, 2010}}

\maketitle
%\vspace{2cm}\noindent
 \begin{abstract}
 This paper examines the existence and region of convergence of Fourier transform of the functions of bicomplex variables with the help of projection on its idempotent components as auxiliary complex planes. Several basic properties of this bicomplex version of Fourier transform are examined. 
 \end{abstract}


\vspace{2cm}
 %\noindent
 %%%%%%%%%%%%%%%%%%%%%%%%%%%%%%%%%%%%%%%%%%%%%%%%%%%%%%%%%%%%%%%%%%%%%%%%%%%%%%%%%%%%%%%%%%%%%%%%%%%%%%%%%%%%%%%%%%%%%%%%%%%%%%%%%%%%%%%%%%%%%%%%%%%

 \section{Introduction}

 The theory of bicomplex numbers is a matter of active research for quite a long time since the seminal work of Segre \cite{segre} in search of a special algebra. The algebra of bicomplex numbers are widely used in the literature as it becomes a viable commutative alternative \cite{spam1,spam2} to the non-commutative skew field of quaternions introduced by Hamilton \cite{ham} (both are four-dimensional and generalization of complex numbers). The commutativity in the former is gained at the cost of the fact that the ring of these numbers contains zero-divisors and so can not form a field \cite{price}. However the novelty of commutativity of bicomplex numbers is that the later can be recognized as the complex numbers with complex coefficients as it's immediate effect and so there are deep similarities between the properties of complex and bicomplex numbers \cite{ol}. Many recent developments have aimed to achieve different algebraic \cite{shp1,shp2,dim,goyal} and geometric \cite{ya1,ya2,cha} properties of bicomplex numbers, the analysis of bicomplex functions \cite{Ga,Moro,Ry,Sri,ronn} and its applications on different branches of physics (such as quantum physics,High energy physics, Bifurcation and chaos etc) \cite{Xu,Krsh,qntm1,roch1,roch2,Mart} to name a few.\\

 In two recent developments \cite{kk,ban} efforts have been done to extend the Laplace transform and its inverse transform in the bicomplex variables from their complex counterpart. In their procedure the idempotent representation of the bicomplex variables plays a vital role. Actually these idempotent components are complex valued and the bicomplex counterpart simply is their combination with idempotent hyperbolic numbers. The Laplace transform of these idempotent complex variables within their regions of convergence are taken and then the bicomplex version of that transform can be obtained directly by combination of them with idempotent hyperbolic numbers. The region of convergence in the later case will be the union of the respective regions of those idempotent complex variables. Bicomplex version of the inversion of Laplace transform is achieved by employing the residual procedure on both the complex planes in connection to idempotent representation.\\

 In the same spirit we take up the study for existence of Fourier transform and the region of convergence in bicomplex variables. The Fourier transform \cite{chand,kaiser} is actually a reversible operation employed to transform signals between the spatial (or time) domain and the frequency domain. Most often in the literature $f$ is a real valued function and its Fourier transform $\hat{f}$ is complex valued where a complex number describes both the amplitude and phase of a corresponding frequency component.\\

 In this paper one of our concern is to extend the Fourier transform in bicomplex variables from its complex version that can capable of transferring  signals from real-valued $(t)$ domain to bicomplex frequency $(\omega)$ domain. The later should have two idempotent complex frequency components $\omega_1 \mbox{ and }\omega_2 $.  \\

 The organization of our paper is as follows: \\
 Section 2 introduces a brief preliminaries of bicomplex numbers. In section 3 we present the existence and region of convergence of bicomplex version of Fourier transform. Some of its basic properties are extended from complex Fourier transform and finally section 4 contains the conclusion.

%%%%%%%%%%%%%%%%%%%%%%%%%%%%%%%%%%%%%%%%%%%%%%%%%%%%%%%%%%%%%%%%%%%%%%%%%%%%%%%%%%%%%%%%%%%%%%%%%%%%%%%%%%%%%%%%%%%%%%%%%%%%%%%%%%%%%%%%%%%%%%%%%%%
\section{Bicomplex numbers}

We start with an unconventional interpretation of the set of complex numbers $\mathbb{C}$ in which its members are found by duplication of the elements of the set of real numbers $\mathbb{R}$ in association with a non-real unit $i$, such that $i^2 =-1$ in the form
\begin{equation}\label{dc}
\mathbb{C}=\{z=x+iy:x,y\in\mathbb{R}\}.
\end{equation}

Now if we repeat our duplication process once on the members of $\mathbb{C}$, for neatness we first denote the imaginary unit $i$ of (\ref{dc}) by $i_1$ resulting \begin{equation}
\mathbb{C}(i_1)=\{z=x+i_1 y:x,y\in\mathbb{R}\}.\nonumber
\end{equation}
If $i_2$ be a new imaginary unit associated with duplication, having the properties
\begin{equation}
{i_2}^2 =-1 ;\quad i_1 i_2 = i_2 i_1 ; \quad ai_2 = i_2 a,\forall a \in \mathbb{R}\nonumber
\end{equation}
we can extend $\mathbb{C}(i_1)$ onto the set of bicomplex numbers
\begin{equation}\label{db}
\mathbb{C}_2 = \{\omega=z_1 +i_2 z_2:z_1 ,z_2 \in \mathbb{C}(i_1)\}
\end{equation}
where an additional structure of commutative multiplication is imbedded.\\

Going back to the real variables, for $z_1 = x_1 +i_1 x_2$ and $z_2 =x_3 +i_1 x_4$, the bicomplex numbers admits of an alternative representation of the form
\begin{equation}
\omega = x_1 + i_1 x_2 +i_2 x_3 +i_1 i_2 x_4 \nonumber
\end{equation}
which is the linear combination of four units: one real unit $1$, two imaginary units $i_1 ,i_2$ and one non-real hyperbolic unit $i_1 i_2 (=i_2 i_1)$ for which
$(i_1 i_2)^2 =1$. In particular if $x_2 =x_3 =0$ one may identify bicomplex numbers with the hyperbolic numbers.\\

However looking onto the algebraic structure of $\mathbb{C}_2$ we can observe that it becomes a commutative ring with unit and $\mathbb{R},\mathbb{C}(i_1)$ are two subrings embedded within it as
\begin{eqnarray}
\mathbb{R} & \equiv & \{z_1 +i_2 z_2 :z_2 =0, z_1 \in \mathbb{R}\}\subset \mathbb{C}_2\nonumber \\
\mathbb{C}(i_1) & \equiv & \{z_1 +i_2 z_2 :z_2 =0, z_1 \in \mathbb{C}(i_1)\}\subset \mathbb{C}_2 .\nonumber
\end{eqnarray}

Interestingly, we may indeed identify the set of complex numbers $\mathbb{C}$ with duplication of reals associated with imaginary unit $i_2$, i.e.
\begin{equation}
\mathbb{C}(i_2)=\{z=x+i_2 y:x,y\in\mathbb{R}\}\nonumber
\end{equation}
 as another possible subring imbedding onto $\mathbb{C}_2$. Both $\mathbb{C}(i_1) \mbox{ and } \mathbb{C}(i_2)$ are isomorphic to $\mathbb{C}$ but are essentially different.\\

 Furthermore for two arbitrary bicomplex numbers $\omega= z_1 +i_2 z_2$ and $\omega'=z'_1+i_2 z'_2$; $z_1,z_2,z'_1,z'_2\in \mathbb{C}(i_1)$ the scalar addition is defined by $$\omega+\omega'=(z_1 +z'_1)+i_2 (z_2 +z'_2)$$ and the scalar multiplication is governed by
 $$\omega . \omega'=(z_1 z'_1 -z_2 z'_2)+i_2 (z_2 z'_1 +z_1 z'_2).$$

 \subsection{Idempotent representation}

 We now introduce two bicomplex numbers
 \begin{equation}\label{idcmp}
 e_1 =\frac{1+i_1 i_2}{2},\quad e_2 =\frac{1-i_1 i_2}{2}
 \end{equation}
  those satisfy
 \begin{equation}
 e_1 +e_2 =1,\quad e_1 . e_2 =e_2 . e_1 =0,\quad {e_1}^2 =e_1 .e_1 =e_1 ,\quad {e_2}^2 =e_2 .e_2 =e_2.\nonumber
 \end{equation}
 The second requirement indicates that $e_1 ,e_2$ are orthogonal while the last two signal them as idempotent. They offer us a unique decomposition of $\mathbb{C}_2$ in the following form:\\
 for any $\omega = z_1 +i_2 z_2 \in \mathbb{C}_2$;
  \begin{equation}\label{decomp}
  z_1 +i_2 z_2 = (z_1 -i_1 z_2 )e_1 + (z_1 +i_1 z_2 )e_2
  \end{equation}
   resulting a pair of mutually complementary projections
   \begin{eqnarray}
   && \mathcal{P}_1 : (z_1 +i_2 z_2 ) \in \mathbb{C}_2 \mapsto (z_1 -i_1 z_2 )\in \mathbb{C}(i_1)\nonumber\\
   && \mathcal{P}_2 : (z_1 +i_2 z_2 ) \in \mathbb{C}_2 \mapsto (z_1 +i_1 z_2 )\in \mathbb{C}(i_1).\nonumber
   \end{eqnarray}

   One may at once verify that ${\mathcal{P}_1 }^2 =\mathcal{P}_1 ,{\mathcal{P}_2 }^2 =\mathcal{P}_2 , \mathcal{P}_1 e_1 +\mathcal{P}_2 e_2 = I$ and for any $\omega_1 ,\omega_2 \in \mathbb{C}_2$;
   \begin{eqnarray}
   \mathcal{P}_k (\omega_1 +\omega_2 ) &=& \mathcal{P}_k (\omega_1 )+\mathcal{P}_k (\omega_2 )\nonumber\\
   \mathcal{P}_k (\omega_1 \omega_2 ) &=& \mathcal{P}_k (\omega_1 )\mathcal{P}_k (\omega_2 ),\quad k=1,2.\nonumber
   \end{eqnarray}

 At this stage, we now mention the auxiliary complex spaces of the space of bicomplex numbers which are
 \begin{eqnarray}
 \mathcal{A}_1 =\{\mathcal{P}_1 (\omega):\omega \in \mathbb{C}_2 \}\nonumber \\
 \mathcal{A}_2 =\{\mathcal{P}_2 (\omega):\omega \in \mathbb{C}_2 \}.\nonumber
 \end{eqnarray}

\subsection{Bicomplex functions}

We start with a bicomplex-valued function $f:\Omega \subset \mathbb{C}_2 \mapsto \mathbb{C}_2$. The derivative of $f$ at a point $\omega_0 \in \Omega$ is defined by
\begin{equation}
f'(\omega_0 )={\lim_{h\rightarrow 0}}\frac{f(\omega_0 +h)-f(\omega_0 )}{h}\nonumber
\end{equation}
provided the limit exists and the domain $\Omega$ is so chosen that $h=h_0 + i_1 h_1  + i_2 h_2  + i_1 i_2 h_3 $ is invertible. (It is direct to prove that $h$ is not invertible only for $h_0 = -h_3 ,h_1 =h_2 \mbox{ or }h_0 = h_3 ,h_1 =-h_2 $).\\

If the bicomplex derivative of $f$ exists at each point of it's domain $\Omega$ then, in similar to complex functions, $f$ will be a bicomplex holomorphic function in $\Omega$. Indeed if $f$ can be expressed as
\begin{equation}
f(\omega)=g_1 (z_1 ,z_2 )+i_2 g_2 (z_1 ,z_2 ), \quad \omega = (z_1 +i_2 z_2) \in \Omega \nonumber
\end{equation}
then $f$ will be holomorphic if and only if $g_1 ,g_2$ are both complex holomorphic in $z_1 ,z_2$ \cite{ronn} and
\begin{equation}
\frac{\partial g_1 }{\partial z_1}= \frac{\partial g_2 }{\partial z_2},\quad \frac{\partial g_1 }{\partial z_2}=-\frac{\partial g_2 }{\partial z_1}.\nonumber
\end{equation}
Moreover $f'(\omega)=\frac{\partial g_1}{\partial z_2}+i_2 \frac{\partial g_2}{\partial z_1}$ and it is invertible only when $\det\left(
                                                                                                                                             \begin{array}{cc}
                                                                                                                                               \frac{\partial g_1}{\partial z_1} &  \frac{\partial
                                                                                                                                               g_1}{\partial z_2}\\
                                                                                                                                               \frac{\partial
                                                                                                                                               g_2}{\partial z_1} & \frac{\partial g_2}{\partial z_2}\\
                                                                                                                                             \end{array}
                                                                                                                                           \right)\neq 0$.
\\

In the following we take up the idempotent representation of bicomplex numbers which is crucial in a deeper understanding of the analysis of holomorphic functions. Any bicomplex holomorphic function $f:\Omega \subset \mathbb{C}_2 \mapsto \mathbb{C}_2$ involving unique idempotent decomposition into two complex- valued functions \cite{ronn} reads as
\begin{equation}
f(\omega)=f_1 (\omega_1 )e_1 + f_2 (\omega_2 )e_2 ,\quad \omega = (\omega_1 e_1 +\omega_2 e_2)\in \Omega. \nonumber
\end{equation}
One may then verify in a straightforward way that
\begin{eqnarray}
\Omega_1 &=& \{\omega_1:\omega \in \Omega\}\subset \mathbb{C}(i_1 )\nonumber\\
\Omega_2 &=& \{\omega_2:\omega \in \Omega\}\subset \mathbb{C}(i_1 )\nonumber
\end{eqnarray}
will be domain of complex-valued functions $f_1$ and $f_2$ respectively. In view of projection operators $\mathcal{P}_1$ and $\mathcal{P}_2$ that can be represented as
\begin{eqnarray}
&&\Omega_1 = \mathcal{P}_1 (\Omega) \quad \Rightarrow \quad f_1 \equiv \mathcal{P}_1 f \nonumber\\
&&\Omega_2 = \mathcal{P}_2 (\Omega) \quad \Rightarrow \quad f_2 \equiv \mathcal{P}_2 f.\nonumber
\end{eqnarray}
\\
Indeed in case of bicomplex-valued holomorphic functions most often the properties of its idempotent complex-valued holomorphic components are just carried over their bicomplex counterpart \cite{Sri}. For example,
  $f(\omega)$ will be convergent in a domain $\Omega$ if and only if $f_1 (\omega_1 ),f_2 (\omega_2 )$ are convergent in their domains $\Omega_1 = \mathcal{P}_1 (\Omega)$ and $\Omega_2 = \mathcal{P}_2 (\Omega)$ respectively.

%%%%%%%%%%%%%%%%%%%%%%%%%%%%%%%%%%%%%%%%%%%%%%%%%%%%%%%%%%%%%%%%%%%%%%%%%%%%%%%%%%%%%%%%%%%%%%%%%%%%%%%%%%%%%%%%%%%%%%%%%%%%%%%%%%%%%%%%%%%%%%%%%
\section{Bicomplex version of Fourier transform}
In this section our aim is to extend the Fourier transform $\mathcal{F}:\mathbb{D} \subset \mathbb{R} \mapsto \mathbb{C}_2$ in bicomplex variables from its complex version and to verify the basic properties in our version those hold good in later case.

\subsection{Conjecture}
Suppose $f(t)$ be a real-valued function that is continuous for $-\infty <t <\infty$ and satisfies the estimates
\begin{eqnarray}\label{ex1}
&& \mid f(t)\mid \leq C_1 \exp(-\alpha t),\quad t\geq 0, \quad \alpha >0 \nonumber\\
&& \mid f(t)\mid \leq C_2 \exp(\beta t),\quad t\leq 0, \quad \beta >0
\end{eqnarray}
which guarantees that $f$ is absolute integrable on the whole real line.\\

Now we start with the complex Fourier transform \cite{sidorov} $\mathcal{F}:\mathbb{D} \subset \mathbb{R}\mapsto \mathbb{C}(i_1)$. The complex Fourier transform of $f(t)$ associted with complex frequency $\omega_1$ is defined by
\begin{equation}
\hat{f}_1 (\omega_1 ) = \mathcal{F}\{f(t)\} = \int_{-\infty}^\infty \exp(i_1 \omega_1 t) f(t)dt,\quad \omega_1 \in \mathbb{C}(i_1)\nonumber
\end{equation}
together with the requirement of $\mid \hat{f}_1 (\omega_1) \mid < \infty$.

Now for $\omega_1 =x+i_1 y$,\begin{eqnarray}
&&\mid \hat{f}_1 (\omega_1 )\mid = \mid \int_{-\infty}^\infty \exp(i_1 \omega_1 t) f(t)dt \mid \leq \int_{-\infty}^\infty \mid \exp(-yt) f(t)\mid dt\nonumber\\
&=& \int_{-\infty}^0 \exp(-yt) \mid f(t) \mid dt+ \int_{0}^\infty \exp(-yt) \mid f(t) \mid dt \nonumber\\
&\leq & C_2 \int_{-\infty}^0 \exp\{(\beta-y)t\} dt +C_1 \int_{0}^\infty \exp\{-(\alpha+y)t\}dt\nonumber\\
&=& C_2 \frac{1}{\beta-y}+C_1 \frac{1}{\alpha+y}\nonumber
\end{eqnarray}
where we use the estimates (\ref{ex1}) and the facts $\mid \exp(i_1 xt) \mid=1,\mid \exp(-yt)\mid=\exp(-yt)$, as $\exp(-yt)>0$.

Then the requirement $\mid \hat{f}_1 (\omega_1) \mid < \infty$ only implies that $-\alpha<y<\beta$. As its consequence $\hat{f}_1 (\omega_1)$ is holomorphic in the strip $$\Omega_1 = \{\omega_1 \in \mathbb{C}(i_1): -\infty < \mbox{ Re }(\omega_1 )<\infty, -\alpha < \mbox{ Im }(\omega_1 )< \beta\}.$$

In similar arguments the complex Fourier transform of $f(t)$ associted with another complex frequency $\omega_2$ will be
\begin{equation}
\hat{f}_2 (\omega_2 ) = \int_{-\infty}^\infty \exp(i_1 \omega_2 t) f(t)dt,\quad \omega_2 \in \mathbb{C}(i_1)\nonumber
\end{equation}
which will be holomorphic in the strip $$\Omega_2 = \{\omega_2 \in \mathbb{C}(i_1): -\infty < \mbox{ Re }(\omega_2 )<\infty, -\alpha < \mbox{ Im }(\omega_2 )< \beta \}.$$

Now employing duplication over these complex functions $\hat{f}_1 (\omega_1 ), \hat{f}_2 (\omega_2 )$ in association with idempotent units $e_1$ and $e_2$ we observe that
\begin{eqnarray}
\hat{f}_1 (\omega_1 )e_1 +\hat{f}_2 (\omega_2 )e_2 &=& \int_{-\infty}^\infty \exp(i_1 \omega_1 t) f(t)dt. e_1 +\int_{-\infty}^\infty \exp(i_1 \omega_2 t) f(t)dt. e_2 \nonumber\\
&=& \int_{-\infty}^\infty \exp(i_1 \{\omega_1 e_1 +\omega_2 e_2\} t)f(t)dt \nonumber \\
&=& \int_{-\infty}^\infty \exp(i_1 \omega t) f(t)dt \nonumber \\
&=& \hat{f}(\omega)\nonumber
\end{eqnarray}
where we use duplication of complex frequencies $\omega_1 ,\omega_2$ to obtain bicomplex frequency $\omega$ as $\omega = \omega_1 e_1 + \omega_2 e_2$.\\

Since $\hat{f}_1 (\omega_1 ),\hat{f}_2 (\omega_2 )$ are complex holomorphic functions in $\Omega_1 ,\Omega_2$ respectively then as it's natural consequence the bicomplex function $\hat{f}(\omega)$ will be holomorphic in the region
\begin{equation}
\Omega = \{\omega \in \mathbb{C}_2 : \omega = \omega_1 e_1 +\omega_2 e_2 , \omega_1 \in \Omega_1 \mbox{ and }\omega_2 \in \Omega_2 \} .\nonumber
\end{equation}

It is worthwhile to mention that the complex-valued holomorphic functions $\hat{f}_1 (\omega_1 ) \mbox{ and } \hat{f}_2 (\omega_2 )$ are both convergent absolutely in $\Omega_1 \mbox{ and } \Omega_2$ respectively. Then it is direct to prove that the region of absolute convergence of $\hat{f}(\omega)$ will be $\Omega$.\\

For better geometrical understanding of the region of convergence of bicomplex Fourier transform it will be advantageous to use the general four-unit representation of bicomplex numbers. In this occasion we take conventional representation of $\omega_1 ,\omega_2 \in \mathbb{C}(i_1)$ as
\begin{equation}\label{rest}
\omega_1 =x_1 +i_1 x_2 , \quad  \omega_2 =y_1 +i_1 y_2 ; \quad x_1 ,x_2 , y_1 ,y_2 \in \mathbb{R}
\end{equation}
where the requirement for $\omega_1 \in \Omega_1$ and $\omega_2 \in \Omega_2$ imply  $ -\infty <x_1 ,y_1 <\infty $\\
 and $-\alpha <x_2 <\beta ;\quad -\alpha <y_2 <\beta$. Using these and (\ref{idcmp}) $\omega$ takes the explicit four-components form
\begin{equation}\label{ex}
\omega =  \frac{x_1 +y_1}{2}+i_1 \frac{x_2 +y_2}{2}+i_2 \frac{y_2 -x_2}{2}+i_1 i_2 \frac{x_1 -y_1}{2}=a_0 +i_1 a_1 +i_2 a_2 +i_1 i_2 a_3
\end{equation}
where $a_0 ,a_1 ,a_2 ,a_3 \in \mathbb{R}$.\\

On the basis of the restrictions on $x_2 \mbox{ and }y_2$ given in (\ref{rest}), the following three possibilities can occur:
\begin{enumerate}
  \item If $x_2 =y_2$, it is trivial to obtain $-\alpha <a_1 <\beta$ and $a_2 =0$,
  \item For $x_2 > y_2$ one may infer $-\alpha -a_2 <a_1 <\beta +a_2$ whereas $-\frac{\alpha + \beta}{2}< a_2 <0$,
  \item If $x_2 < y_2$ then in similar to previous possibility we obtain $-\alpha +a_2 < a_1 <\beta - a_2$ and $0 < a_2 < \frac{\alpha +\beta}{2}$,
\end{enumerate}
whereas $-\infty <a_0 ,a_3 <\infty$ in all three cases.\\

Considering all of these results we conclude that $$-\infty <a_0 ,a_3 <\infty, -\alpha +\mid a_2 \mid <a_1 <\beta -\mid a_2 \mid \mbox{ and }0\leq \mid a_2 \mid <\frac{\alpha + \beta}{2}$$
and hence the region of convergence of $\hat{f}(\omega)$ (See  fig-1 in appendix for $a_1 - a_2$ plane section of the region) can be identified as
\begin{equation}\label{roc}
\Omega = \{\omega =a_0 +i_1 a_1 +i_2 a_2 +i_1 i_2 a_3 \in \mathbb{C}_2 : -\infty <a_0 ,a_3 <\infty; -\alpha +\mid a_2 \mid <a_1 <\beta -\mid a_2 \mid; 0\leq \mid a_2 \mid <\frac{\alpha + \beta}{2}\}.
\end{equation}

\begin{itemize}
\item Conversely, the existence  of bicomplex Fourier transform $\hat{f}(\omega)$ can be obtained in the following way:\\
 If  $\omega = a_0 +i_1 a_1 +i_2 a_2 +i_1 i_2 a_3 \in \Omega$ ; $-\infty <a_0 ,a_3 <\infty, -\alpha +\mid a_2 \mid  <a_1 <\beta -\mid a_2 \mid, \mbox{ and } 0\leq \mid a_2 \mid <\frac{\alpha + \beta}{2}$.
 Now expressing $\omega$ in idempotent components as
 $$\omega= a_0 +i_1 a_1 +i_2 a_2 +i_1 i_2 a_3 = [(a_0 +a_3 )+i_1 (a_1 -a_2 )]e_1 + [(a_0 -a_3 )+i_1 (a_1 +a_2 )]e_2 =\omega_1 e_1 +\omega_2 e_2$$
 we obtain
 \begin{enumerate}
   \item $a_2 =0$ and $-\alpha <a_1 <\beta$ trivially leads $-\alpha< a_1 -a_2 <\beta$ and $-\alpha<a_1 +a_2 <\beta$,
   \item when $a_2 <0$, from the first inequality of $-\alpha-a_2 <a_1 <\beta +a_2$ we can get $-\alpha<a_1 +a_2$ whereas the last inequality gives $a_1 -a_2 <\beta$. Following $a_2 <0$ these results can be interpreted as $-\alpha<a_1 +a_2 <a_1 -a_2$ and $a_1 +a_2 <a_1 -a_2 <\beta$ which in together combined into $-\alpha <a_1 +a_2 <a_1 -a_2 <\beta$,
   \item when $a_2 >0$, from the first inequality of $-\alpha +a_2 <a_1 <\beta -a_2$ we can get $-\alpha<a_1 -a_2$ whereas the last inequality gives $a_1 +a_2 <\beta$. Following $a_2 >0$ these results can be interpreted as $-\alpha<a_1 -a_2 <a_1 +a_2$ and $a_1 -a_2 <a_1 +a_2 <\beta$ which in together combined into $-\alpha <a_1 -a_2 <a_1 +a_2 <\beta$,
 \end{enumerate}
 Hence the result.
\end{itemize}

Now we are ready to define the Fourier transform for bicomplex variable.

\subsection{Definition}
Let $f(t)$ be a real-valued continuous function in $(-\infty,\infty)$ that satisfies the estimates  (\ref{ex1}).
The Fourier transform of $f(t)$ can be defined as
\begin{equation}\label{ft}
\hat{f}(\omega ) = \mathcal{F}\{f(t)\} = \int_{-\infty}^\infty \exp(i_1 \omega t) f(t)dt,\quad \omega \in \mathbb{C}_2 .
\end{equation}
The Fourier transform $\hat{f}(\omega)$ exists and holomorphic for all $\omega \in \Omega$ where $\Omega$ (given in (\ref{roc})) is the region of absolute convergence of $\hat{f}$.

\subsection{Existence of Fourier transform}
\textbf{\textit{Theorem:}} If $f(t)$ be a real valued function and is continuous for $-\infty<t<\infty$ satisfying estimates (\ref{ex1}) then $\hat{f}(\omega)$ (defined in (\ref{ft})) exists in the region (\ref{roc}).\\\\

\textbf{\textit{Proof:}}\begin{eqnarray}
\hat{f}(\omega)&=& \int_{-\infty}^\infty \exp(i_1 \omega t) f(t)dt\nonumber\\
&=& \int_{-\infty}^\infty \exp(i_1 \omega_1 t) f(t)dt. e_1 + \int_{-\infty}^\infty \exp(i_1 \omega_2 t) f(t)dt .e_2\nonumber
 \end{eqnarray}
 Both the integrals exist when $-\alpha<\mbox{ Im }(\omega_1 =x_1 +i_1 x_2 )<\beta$ and $-\alpha<\mbox{ Im }(\omega_2 =y_1 +i_1 y_2 )<\beta$. So $\hat{f}(\omega)$ exists for $\omega = \omega_1 e_1 +\omega_2 e_2 =a_0 +i_1 a_1 +i_2 a_2 +i_1 i_2 a_3$ where $$-\infty <a_0 ,a_3 <\infty, -\alpha +\mid a_2 \mid  <a_1 <\beta -\mid a_2 \mid, \mbox{ and } 0\leq \mid a_2 \mid <\frac{\alpha + \beta}{2}.$$

 \subsection{Uniqueness of Fourier transform}
\textbf{\textit{Theorem:}} If $f(t)$ and $g(t)$ have Fourier transforms $\hat{f}(\omega)$ and $\hat{g}(\omega)$ respectively and \\ $\hat{f}(\omega)=\hat{g}(\omega)$, then $f(t)=g(t)$. \\\\

\textbf{\textit{Proof:}} Let $\hat{f}(\omega)=\hat{f}_1 (\omega_1 )e_1 +\hat{f}_2 (\omega_2 )e_2$ and $\hat{g}(\omega)=\hat{g}_1 (\omega_1 )e_1 +\hat{g}_2 (\omega_2 )e_2$ in their idempotent representations. Now $\hat{f}(\omega)=\hat{g}(\omega)$ is possible if and only if
\begin{eqnarray}
&& \hat{f}_1 (\omega_1 )=\hat{g}_1 (\omega_1 )\mbox{ and }\hat{f}_2 (\omega_2 )=\hat{g}_2 (\omega_2 )\nonumber\\
& \Rightarrow & \int_{-\infty}^\infty \exp(i_1 \omega_1 t) f(t)dt = \int_{-\infty}^\infty \exp(i_1 \omega_1 t) g(t)dt \nonumber\\
&& \mbox{ and } \int_{-\infty}^\infty \exp(i_1 \omega_2 t) f(t)dt = \int_{-\infty}^\infty \exp(i_1 \omega_2 t) g(t)dt\nonumber\\
& \Rightarrow & f(t)=g(t).\nonumber
\end{eqnarray}

\subsection{Basic properties of Fourier transform}
\begin{enumerate}
  \item \textbf{Linearity property}\\
   \textbf{\textit{Theorem:}} If the Fourier transforms of $f(t)$ and $g(t)$ are $\hat{f}(\omega)$ and $\hat{g}(\omega)$ respectively and $a$ and $b$ are constants then $\mathcal{F}(af+bg)=a\hat{f}+b\hat{g}$.

   \textbf{\textit{Proof:}} Let $\hat{f}$ and $\hat{g}$ are both defined for $\omega(=\omega_1 e_1 +\omega_2 e_2 ) \in \Omega,\quad \omega_1 \in \Omega_1, \omega_2 \in \Omega_2$, then
   \begin{eqnarray}
   && \mathcal{F}\left\{a f(t)+b g(t)\right\}=\int_{-\infty}^\infty \exp(i_1 \omega t)\left\{a f(t)+b g(t)\right\}dt\nonumber\\
   &=& \int_{-\infty}^\infty \exp(i_1 \{\omega_1 e_1 +\omega_2 e_2\} t)\left\{a f(t)+b g(t)\right\}dt\nonumber\\
   &=& \int_{-\infty}^\infty \exp(i_1 \omega_1 t)\left\{a f(t)+b g(t)\right\}dt.e_1 +\int_{-\infty}^\infty \exp(i_1 \omega_2 t)\left\{a f(t)+b g(t)\right\}dt.e_2 \nonumber\\
   &=& a\left[\int_{-\infty}^\infty \exp(i_1 \omega_1 t)f(t)dt.e_1 +\int_{-\infty}^\infty \exp(i_1 \omega_2 t)f(t)dt.e_2\right]\nonumber\\
   && +b\left[\int_{-\infty}^\infty \exp(i_1 \omega_1 t)g(t)dt.e_1 +\int_{-\infty}^\infty \exp(i_1 \omega_2 t)g(t)dt.e_2 \right]\nonumber\\
   &=& a\int_{-\infty}^\infty \exp(i_1 \{\omega_1 e_1 +\omega_2 e_2\} t)f(t)dt +b\int_{-\infty}^\infty \exp(i_1 \{\omega_1 e_1 +\omega_2 e_2\} t)g(t)dt\nonumber\\
   &=& a\int_{-\infty}^\infty \exp(i_1 \omega t)f(t)dt +b\int_{-\infty}^\infty \exp(i_1 \omega t)g(t)dt = a\hat{f}(\omega)+b\hat{g}(\omega).\nonumber
   \end{eqnarray}

\item \textbf{Shifting property}\\
   \textbf{\textit{Theorem:}} If $\hat{f}(\omega)$ is the Fourier transforms of $f(t)$ then $\mathcal{F}\left\{f(t-a)\right\}=\exp(i_1 \omega a)\hat{f}(\omega)$.

   \textbf{\textit{Proof:}} By definition $\mathcal{F}\left\{f(t-a)\right\}=\int_{-\infty}^\infty \exp(i_1 \omega t)f(t-a)dt$.\\
   Now for $t=a+u$, the integral in right hand side is equal to \\ $\int_{-\infty}^\infty \exp\{i_1 \omega (a+u)\}f(u)du=\exp(i_1 \omega a)\int_{-\infty}^\infty \exp(i_1 \omega u)f(u)du=\exp(i_1 \omega a)\hat{f}(\omega)$.

   \item \textbf{Scaling property}\\
   \textbf{\textit{Theorem:}} If $\hat{f}(\omega)$ is the Fourier transforms of $f(t)$ then $\mathcal{F}\{f(at)\}=\frac{1}{\mid a \mid}\hat{f}(\frac{\omega}{a})$ where $a\neq 0$.

   \textbf{\textit{Proof:}} If $a>0$ then $ \mathcal{F}\{f(at)\}= \int_{-\infty}^\infty \exp(i_1 \omega t)f(at)dt =\frac{1}{a}\int_{-\infty}^\infty \exp(i_1 \frac{\omega}{a} u)f(u)du = \frac{1}{a}\hat{f}(\frac{\omega}{a})$ where we take $at=u$.

   If $a<0$ then for $a=-b :b>0$ we have $\mathcal{F}\{f(at)\}=\int_{-\infty}^\infty \exp(i_1 \omega t)f(at)dt=\int_{-\infty}^\infty \exp(i_1 \omega t)f(-bt)dt$. Now taking $bt=-u$ the integral is $-\frac{1}{b}\int_{\infty}^{-\infty} \exp(i_1 \frac{\omega}{-b} u)f(u)du =\frac{1}{b}\int_{-\infty}^\infty \exp(i_1 \frac{\omega}{-b} u)f(u)du =\frac{1}{-a}\int_{-\infty}^\infty \exp(i_1 \frac{\omega}{a} u)f(u)du = \frac{1}{-a}\hat{f}(\frac{\omega}{a})$.

   From above these results we conclude $\mathcal{F}\{f(at)\}=\frac{1}{\mid a \mid}\hat{f}(\frac{\omega}{a})$.

 \item \textbf{Convolution theorem}\\
   \textbf{\textit{Theorem:}} The Fourier transform of the convolution of two functions $f(t)$ and $g(t)$,$-\infty<t<\infty$ is the product of their Fourier transforms, respectively $\hat{f}(\omega)$ and $\hat{g}(\omega)$ \textit{i.e.}
   $$\mathcal{F}\left\{f(t)*g(t)\right\}=\mathcal{F}\left\{\int_{-\infty}^\infty f(u)g(t-u)du\right\}=\hat{f}(\omega)\hat{g}(\omega).$$

   \textbf{\textit{Proof:}} By definition
   \begin{eqnarray}
   && \mathcal{F} \left\{f(t)*g(t)\right\}=\mathcal{F}\left\{\int_{-\infty}^\infty f(u)g(t-u)du\right\}=\int_{-\infty}^\infty \exp(i_1 \omega t)\left\{\int_{-\infty}^\infty f(u)g(t-u)du\right\}dt\nonumber\\
   &=&\int_{-\infty}^\infty \exp(i_1 \omega_1 t)\left\{\int_{-\infty}^\infty f(u)g(t-u)du\right\}dt. e_1 +\int_{-\infty}^\infty \exp(i_1 \omega_2 t)\left\{\int_{-\infty}^\infty f(u)g(t-u)du\right\}dt. e_2\nonumber\\
   &=&\int_{-\infty}^\infty f(u)\left\{\int_{-\infty}^\infty \exp(i_1 \omega_1 t)g(t-u)dt\right\}du. e_1 +\int_{-\infty}^\infty f(u)\left\{\int_{-\infty}^\infty \exp(i_1 \omega_2 t)g(t-u)dt\right\}du. e_2\nonumber\\
   && \mbox{ using method for changing order of integrals in complex analysis \cite{mat}}\nonumber\\
   &=&\int_{-\infty}^\infty f(u)\left\{\int_{-\infty}^\infty \exp(i_1 \omega t)g(t-u)dt\right\}du\nonumber\\
   &=& \int_{-\infty}^\infty f(u)\exp(i_1 \omega u)\hat{g}(\omega)du, \mbox{ using shifting property (see property 2) }\nonumber\\
   &=&  \left\{\int_{-\infty}^\infty f(u)\exp(i_1 \omega u)du\right\}\hat{g}(\omega) = \hat{f}(\omega)\hat{g}(\omega).\nonumber
   \end{eqnarray}

\item \textbf{\textit{Theorem:}} If $f(t)$ and $t^r f(t),r=1,2,....,n$ are all integrable in $-\infty<t<\infty$ then
$$\mathcal{F}\left\{t^n f(t)\right\}=(-i_1)^n \frac{d^n}{d\omega^n}\{\hat{f}(\omega)\}$$
where $\hat{f}(\omega)$ is the Fourier transform of $f(t)$.

   \textbf{\textit{Proof:}} We will prove this theorem by using the method of mathematical induction and differentiation under integral sign. For $n=1$,
   \begin{eqnarray}
   && \frac{d}{d\omega}\hat{f}(\omega) = \frac{\partial}{\partial\omega_1}\hat{f_1}(\omega_1)e_1 +\frac{\partial}{\partial\omega_2}\hat{f_2}(\omega_2)e_2 \nonumber\\
   &=& \frac{\partial}{\partial\omega_1}\int_{-\infty}^\infty \exp(i_1 \omega_1 t)f(t)dt.e_1 + \frac{\partial}{\partial\omega_2}\int_{-\infty}^\infty \exp(i_1 \omega_2 t)f(t)dt.e_2\nonumber\\
   &=& \int_{-\infty}^\infty \frac{\partial}{\partial\omega_1}\{\exp(i_1 \omega_1 t)f(t)\}dt.e_1 + \int_{-\infty}^\infty \frac{\partial}{\partial\omega_2}\{\exp(i_1 \omega_2 t)f(t)\}dt.e_2  \nonumber\\
   && \mbox{using Leibnitz rule in complex analysis \cite{mat} }\nonumber\\
   &=& i_1 \int_{-\infty}^\infty t f(t) \exp(i_1 \omega_1 t)dt.e_1 +i_1 \int_{-\infty}^\infty t f(t) \exp(i_1 \omega_2 t)dt.e_2\nonumber\\
    &=& i_1 \int_{-\infty}^\infty t f(t)\{\exp(i_1 \omega_1 t)e_1 + \exp(i_1 \omega_2 t)e_2\}dt \nonumber\\
    &=&  i_1 \int_{-\infty}^\infty t f(t) \exp(i_1 \omega t)dt=i_1 \mathcal{F}\{t f(t)\}\nonumber\\
   &\Rightarrow & \mathcal{F}\{t f(t)\}=-i_1 \frac{d}{d\omega}\hat{f}(\omega).\nonumber
   \end{eqnarray}
Now for $n=2$, in similar to the case for $n=1$,
\begin{eqnarray}
&& \frac{d^2}{d\omega^2}\hat{f}(\omega)=\frac{d}{d\omega}\left[\frac{d}{d\omega}\hat{f}(\omega)\right]\nonumber\\
&=& i_1 \frac{d}{d\omega}\left[\int_{-\infty}^\infty t f(t)\exp(i_1 \omega t)dt\right]\nonumber\\
&=& i_1 \frac{\partial}{\partial\omega_1}\int_{-\infty}^\infty \exp(i_1 \omega_1 t)t f(t)dt.e_1 + i_1 \frac{\partial}{\partial\omega_2}\int_{-\infty}^\infty \exp(i_1 \omega_2 t)t f(t)dt.e_2\nonumber\\
   &=& i_1 \int_{-\infty}^\infty \frac{\partial}{\partial\omega_1}\{\exp(i_1 \omega_1 t) t f(t)\}dt.e_1 + i_1 \int_{-\infty}^\infty \frac{\partial}{\partial\omega_2}\{\exp(i_1 \omega_2 t) t f(t)\}dt.e_2  \nonumber\\
   && \mbox{using Leibnitz rule \cite{mat}}\nonumber\\
   &=& - \int_{-\infty}^\infty t^2 f(t) \exp(i_1 \omega_1 t)dt.e_1 - \int_{-\infty}^\infty t^2 f(t) \exp(i_1 \omega_2 t)dt.e_2\nonumber\\
    &=& - \int_{-\infty}^\infty t^2 f(t) \{\exp(i_1 \omega_1 t)e_1 + \exp(i_1 \omega_2 t)e_2\}dt \nonumber\\
    &=& - \int_{-\infty}^\infty t^2 f(t) \exp(i_1 \omega t)dt= -\mathcal{F}\{t^2 f(t)\}\nonumber\\
&\Rightarrow & \mathcal{F}\left\{t^2 f(t)\right\}=-\frac{d^2}{d\omega^2}\hat{f}(\omega)=(-i_1 )^2 \frac{d^2}{d\omega^2}\hat{f}(\omega).\nonumber
\end{eqnarray}
Proceeding in this way we obtain
$$ \mathcal{F}\left\{t^n f(t)\right\}=(-i_1 )^n \frac{d^n}{d\omega^n}\hat{f}(\omega).$$

\item \textbf{\textit{Theorem:}} If $f(t)$ and $f^{(r)}(t),r=1,2,....,n$ are piecewise smooth and tend to 0 as $\mid t \mid \rightarrow \infty$, and $f$ with its derivatives of order up to $n$ are integrable in $-\infty<t<\infty$ then
$$\mathcal{F}\{f^{(n)}(t)\}=(-i_1 \omega)^n \hat{f}(\omega)\}$$
where $\hat{f}(\omega)$ is the Fourier transform of $f(t)$ and $f^{(r)}(t)=\frac{d^r}{dt^r}f(t)$.

   \textbf{\textit{Proof:}} We will prove it also using method of induction. For $n=1$,
   \begin{eqnarray}
   && \mathcal{F}\{f'(t)\}=\int_{-\infty}^\infty \exp(i_1 \omega t)f'(t)dt\nonumber\\
   &=& \left[f(t)\exp(i_1 \omega t)\right]_{-\infty}^\infty -i_1 \omega \int_{-\infty}^\infty \exp(i_1 \omega t) f(t)dt\nonumber\\
   &=& 0-i_1 \omega \hat{f}(\omega)=-i_1 \omega \hat{f}(\omega).\nonumber
   \end{eqnarray}
   Similarly for $n=2$,
   \begin{eqnarray}
   && \mathcal{F}\{f''(t)\}=\int_{-\infty}^\infty \exp(i_1 \omega t)f''(t)dt\nonumber\\
   &=& \left[f'(t)\exp(i_1 \omega t)\right]_{-\infty}^\infty -i_1 \omega \int_{-\infty}^\infty \exp(i_1 \omega t) f'(t)dt\nonumber\\
   &=& 0-i_1 \omega \mathcal{F}\{f'(t)\}=(-i_1 \omega)^2 \hat{f}(\omega).\nonumber
   \end{eqnarray}
   Proceeding with similar arguments we can get $\mathcal{F}\{f^{(n)}(t)\}=(-i_1 \omega)^n \hat{f}(\omega)\}$.

\end{enumerate}

\begin{itemize}
\item \textbf{Corollary}\\
      If $f(t)$ is finite, \textit{i.e.} $f(t)=0 \quad \mid t \mid >T$ and continuous inside $\mid t \mid \leq T$, then its complex Fourier transform is an entire function. As it's consequence $\hat{f}_1 (\omega_1)=\int_{-T}^T \exp(i_1 \omega_1 t)f(t)dt$ and $\hat{f}_2 (\omega_2)=\int_{-T}^T \exp(i_1 \omega_2 t)f(t)dt$. So the bicomplex Fourier transform $\hat{f}(\omega)$ exists and converges absolutely within the whole $\mathbb{C}_2$.
\end{itemize}

\subsection{Examples}
\begin{enumerate}
  \item If \begin{equation*}
f(t)=\exp(-a\mid t \mid),\quad a>0
\end{equation*}
then it satisfies estimates (\ref{ex1}) for $\alpha = \beta =a$ and its complex Fourier transforms are
\begin{equation*}
\hat{f}_1 (\omega_1)=\frac{2a}{a^2 +{\omega_1}^2},\quad \hat{f}_2 (\omega_2)=\frac{2a}{a^2 +{\omega_2}^2}.
\end{equation*}
Both $\hat{f}_1 ,\hat{f}_2$ are holomorphic in the strip $-a<\mbox{ Im }(\omega_1),\mbox{ Im }(\omega_2)<a$.
Then the bicomplex Fourier transform will be
\begin{equation*}
\hat{f}(\omega)=\frac{2a}{a^2 +\omega^2}
\end{equation*}
 with region of convergence $$\Omega=\{\omega=a_0 +i_1 a_1 +i_2 a_2 +i_1 i_2 a_3 \in \mathbb{C}_2 :\quad 0\leq \mid a_2 \mid <a,\quad -a+\mid a_2 \mid <a_1 <
 a-\mid a_2 \mid \}.$$

 \item If $$f(t)=\left\{
              \begin{array}{ll}
                \exp(-t),\qquad{t>0;} \\
               0, \qquad{t\leq 0}
              \end{array}
            \right.$$
then $\alpha=1$ but $\beta$ be any positive number.
Here $$\hat{f}(\omega)=\frac{1}{1-i_1 \omega}$$ and its region of convergence is $$\Omega=\{\omega=a_0 +i_1 a_1 +i_2 a_2 +i_1 i_2 a_3 \in \mathbb{C}_2 : 0\leq \mid a_2 \mid <\frac{1+\beta}{2},a_1 >-1,\quad \mbox{for any positive }\beta\}.$$

\item If $f(t)=\exp(-\frac{t^2}{2})$ then $\hat{f}(\omega)=\sqrt{2\pi}\exp(-\frac{\omega^2}{2})$ and its region of convergence is $$\Omega=\{\omega=a_0 +i_1 a_1 +i_2 a_2 +i_1 i_2 a_3 \in \mathbb{C}_2 : -\infty < a_1, a_2 <\infty \}.$$

\item If $$f(t)=\left\{
                 \begin{array}{ll}
                   1, \quad{\mid t \mid \leq a;} \\
                   0, \quad{\mid t \mid >a}
                 \end{array}
               \right.
$$
then its complex Fourier transforms in both $\omega_1$ and $\omega_2$ planes are entire functions. Then using the corollary we obtain
$\hat{f}_1 (\omega_1)=\frac{2}{\omega_1}\sin(a\omega_1)$ and $\hat{f}_2 (\omega_2)=\frac{2}{\omega_2}\sin(a\omega_2)$, where the singularity at $\omega=0$ is removable. In this case the bicomplex Fourier transform will be $\hat{f}(\omega)=\frac{2}{\omega}\sin(a\omega)$ and it's region of convergence is $\mathbb{C}_2$.

\item If $$f(t)=\left\{
                 \begin{array}{ll}
                   0, \quad{t<0;} \\
                   \exp (-\frac{t}{T})\sin (\omega_0 t), \quad{t\geq 0, \quad T,\omega_0 >0}
                 \end{array}
               \right.$$
which might represent the displacement of a damped harmonic oscillator. Here from the estimates (\ref{ex1}) we have $\alpha=\frac{1}{T}$. Then complex Fourier transform in $\omega_1$ (similar for $\omega_2$) plane is given by $$\hat{f}_1 (\omega_1)=\frac{1}{2}\left[\frac{1}{\omega_1 -\omega_0 +\frac{i_1}{T}}-\frac{1}{\omega_1 +\omega_0 +\frac{i_1}{T}}\right]$$ which is holomorphic in in the infinite strip $\mbox{Im }(\omega_1)>-\frac{1}{T}$ except $\mbox{Re }(\omega_1) \neq \pm \omega_0$. In this problem the bicomplex Fourier transform will be
$$\hat{f} (\omega)=\frac{1}{2}\left[\frac{1}{\omega +\omega_0 +\frac{i_1}{T}}-\frac{1}{\omega -\omega_0 +\frac{i_1}{T}}\right]$$
with region of convergence
$$\Omega=\{\omega =a_0 +i_1 a_1  +i_2 a_2 +i_1 i_2 a_3\in \mathbb{C}_2 : a_0 \neq 0,\pm \omega_0; a_3 \neq 0,\pm \omega_0 ; a_1 >-\frac{1}{T}\}$$ where $a_2 =\frac{\mbox{ Im }(\omega_2)-\mbox{ Im }(\omega_1)}{2}: \mbox{ Im }(\omega_1),\mbox{ Im }(\omega_2)>-\frac{1}{T}$.
\end{enumerate}

\section{Conclusion}

In this paper we have exploited the bicomplex version of Fourier transform method and the condition of absolute convergence of the transformed bicomplex-valued function. We examine the usual properties of Fourier transform in the ring of bicomplex numbers. Finally let us point out that in our observation the concepts introduced or  results obtained in this work are the generalization of the corresponding concepts or results in complex analysis.

 %%%%%%%%%%%%%%%%%%%%%%%%%%%%%%%%%%%%%%%%%%%%%%%%%%%%%%%%%%%%%%%%%%%%%%%%%%%%%%%%%%%%%%%%%%%%%%%%%%%%%%%%%%%%%%%%%%%%%%%%%%%%%%%%%%%%%%%%%%%%%%%%%%%%%%%%%%%
 
%%%%%%%%%%%%%%%%%%%%%%%%%%%%%%%%%%%%%%%%%%%%%%%%%%%%%%%%%%%%%%%%%%%%%%%%%%%%%%%%%%%%%%%%%%%%%%%%%%%%%%%%%%%%%%%%%%%%%%%%%%%%%%%%%%%%%%%%%%%%%%%%%%%%%%%%%%

\begin{thebibliography}{99}
 \bibitem{segre} \textbf{C.Segre} Math.Ann.:40,1892,pp:413.
 \bibitem{spam1} \textbf{N.Spampinato} Atti Reale Accad. Naz. Lincei,Rend.:22,1935,pp:38.
 \bibitem{spam2} \textbf{N.Spampinato} Ann.Mat.Pura.Appl.:14,1936,pp:305.
 \bibitem{ham} \textbf{W.R.Hamilton} \textit{Lectures on quaternion} Dublin:Hodges and Smith: 1853.
 \bibitem{price} \textbf{G.B.Price} Marcel,Dekkar: 1991.
 \bibitem{ol} \textbf{S.Olariu} Norh-Holland Mathematics Studies,Elsevier: 190,2002,pp:269.
 \bibitem{shp1} \textbf{G.Shpilker} Doklady AN SSSR: 282,1985,pp:1090.
 \bibitem{shp2} \textbf{G.Shpilker} Doklady AN SSSR: 293,1987,pp:578.
 \bibitem{dim} \textbf{S.Dimiev,R.Lazov,S.Slavova} Topics in Contemporary Differential Geometry,Complex Analysis and Mathematical physics: 2006,pp:50-56.
 \bibitem{ya1} \textbf{I.M.Yaglom} Academic Press, N.Y.: 1968.
 \bibitem{ya2} \textbf{I.M.Yaglom} Springer, N.Y.: 1979.
 \bibitem{cha} \textbf{K.S.Charak,D.Rochon,N.Sharma} Fractals: 17,2009.
 \bibitem{goyal} \textbf{R.Goyal} Tokyo Journal of Mathematics:30,2007.
 \bibitem{Ga} \textbf{S.Gal} Nova Science Publishers: 2002.
 \bibitem{Moro} \textbf{A.Motter,M.Rosa} Adv. Appl. Clifford Algebra: 8,1998,pp:109.
 \bibitem{Ry} \textbf{J.Ryan} Complex variables,Theory and Applications: 1,1982,pp:119.
 \bibitem{Sri} \textbf{R.K.Srivastava} Proc.Soc.of Special Functions and their applications (SSFA): 2005,pp:55.
 \bibitem{ronn} \textbf{S.R\"{o}nn} arXiv:math/0101200[math.CV],2001.
 \bibitem{Xu} \textbf{Y.Xuegang } Adv.Appl.Clifford Algebra: 9,1998,pp:109.
 \bibitem{Krsh} \textbf{V.Kravchenko,M.Shapiro} Pitman Research Notes in Math.,Addison-Wesley-Longman: 351,1996.
 \bibitem{qntm1} \textbf{D.Dart,D.Haag,H.Cartarins,J.Main,G.Wunner} arXiv:1306.3871[quant-ph],2013.
 \bibitem{roch1} \textbf{D.Rochon,S.Tremblay} Adv.Appl.Clifford Algebra: 14,2004,pp:231.
 \bibitem{roch2} \textbf{D.Rochon,S.Tremblay} Adv.Appl.Clifford Algebra: 16,2006,pp:135.
 \bibitem{Mart} \textbf{E.Martineau,D.Rochon} Int. J. Bifurcation Chaos: 15,2005.
 \bibitem{kk} \textbf{A.Kumar,P.Kumar}  International Journal of Engineering and Technology: 3,2011,pp:225.
 \bibitem{ban}\textbf{A.Banerjee,S.K.Datta,A.Hoque} Mathematical Inverse Problems:1,2014.
  \bibitem{chand} \textbf{S.Bochner,K.Chandrasekharan} Princeton University Press: 1949.
 \bibitem{kaiser} \textbf{G.Kaiser} Birkhauser: 1994.
 \bibitem{sidorov} \textbf{Y.V.Sidorov,M.V.Fedoryuk,M.I.Shabunin} Mir Publishers,Moscow: 1985.
 \bibitem{mat}\textbf{J.H.Mathews,R.W.Howell} Narosa Publication : 2006.
 \end{thebibliography}
 \end{document}